\documentclass{amsart}[12pt]
\usepackage[]{amsmath, amsthm, amsfonts, verbatim, amssymb}
\usepackage{tikz}
\usetikzlibrary{matrix,arrows}
\usepackage{url}
\usepackage[all,cmtip]{xy}

\usepackage{xcolor}

\usepackage{graphicx}
\usepackage{pinlabel}

\def\cH{\mathcal H}

\def\cU{{\mathcal U}}
\def\cV{{\mathcal V}}

\def\bZ{{\mathbb Z}}
\def\bQ{{\mathbb Q}}
\def\bC{{\mathbb C}}

\begin{document}
\newtheorem {theo}{Theorem}
\newtheorem {coro}{Corollary}
\newtheorem {lemm}{Lemma}
\newtheorem {rem}{Remark}
\newtheorem {defi}{Definition}
\newtheorem {ques}{Question}
\newtheorem {prop}{Proposition}
\def\spb{\smallpagebreak}
\def\mpb{\vskip 0.5truecm}
\def\bpb{\vskip 1truecm}
\def\wtM{\widetilde M}
\def\tM{\widetilde M}
\def\wtN{\widetilde N}
\def\tN{\widetilde N}
\def\tR{\widetilde R}
\def\tC{\widetilde C}
\def\tX{\widetilde X}
\def\tY{\widetilde Y}
\def\tE{\widetilde E}
\def\tH{\widetilde H}
\def\tL{\widetilde L}
\def\tQ{\widetilde Q}
\def\tS{\widetilde S}
\def\tc{\widetilde c}
\def\talpha{\widetilde\alpha}
\def\ti{\widetilde \iota}
\def\hM{\hat M}
\def\hq{\hat q}
\def\hR{\hat R}
\def\bs{\bigskip}
\def\ms{\medskip}
\def\ni{\noindent}
\def\td{\nabla}
\def\pd{\partial}
\def\hol{$\text{hol}\,$}
\def\Log{\mbox{Log}}
\def\bfQ{{\bf Q}}
\def\Todd{\mbox{Todd}}
\def\top{\mbox{top}}
\def\Pic{\mbox{Pic}}
\def\bP{{\mathbb P}}
\def\dxi{d x^i}
\def\dxj{d x^j}
\def\dyi{d y^i}
\def\dyj{d y^j}
\def\dzi{d z^I}
\def\dzj{d z^J}
\def\ozi{d{\overline z}^I}
\def\ozj{d{\overline z}^J}
\def\oz1{d{\overline z}^1}
\def\oz2{d{\overline z}^2}
\def\oz3{d{\overline z}^3}
\def\sI{\sqrt{-1}}
\def\hol{$\text{hol}\,$}
\def\ok{\overline k}
\def\ol{\overline l}
\def\oJ{\overline J}
\def\oT{\overline T}
\def\oS{\overline S}
\def\oV{\overline V}
\def\oW{\overline W}
\def\oY{\overline Y}
\def\oL{\overline L}
\def\oI{\overline I}
\def\oK{\overline K}
\def\oL{\overline L}
\def\oj{\overline j}
\def\oi{\overline i}
\def\ok{\overline k}
\def\oz{\overline z}
\def\om{\overline mu}
\def\on{\overline nu}
\def\oa{\overline \alpha}
\def\ob{\overline \beta}
\def\oGamma{\overline \Gamma}
\def\of{\overline f}
\def\oN{\overline N}
\def\og{\overline \gamma}
\def\ogamma{\overline \gamma}
\def\odelta{\overline \delta}
\def\otheta{\overline \theta}
\def\ophi{\overline \phi}
\def\opd{\overline \partial}
\def\oA{\overline A} 
\def\oB{\overline B}
\def\oC{\overline C}
\def\oD{\overline D}
\def\oIq1{\oI_1\cdots\oI_{q-1}}
\def\oIq2{\oI_1\cdots\oI_{q-2}}
\def\op{\overline \partial}
\def\ua{{\underline {a}}}
\def\us{{\underline {\sigma}}}
\def\Chow{{\mbox{Chow}}}
\def\vol{{\mbox{vol}}}
\def\dim{{\mbox{dim}}}
\def\rank{{\mbox{rank}}}
\def\diag{{\mbox{diag}}}
\def\tor{\mbox{tor}}
\def\supp{\mbox supp}
\def\bp{{\bf p}}
\def\bk{{\bf k}}
\def\a{{\alpha}}
\def\tchi{\widetilde{\chi}}
\def\ta{\widetilde{\alpha}}
\def\ovarphi{\overline \varphi}
\def\ocH{\overline{\cH}}
\def\tV{\widetilde{V}}
\def\tf{\widetilde{f}}
\def\th{\widetilde{h}}
\def\tT{\widetilde T}
\def\hG{\widehat{G}}
\def\hS{\widehat{S}}
\def\hD{\widehat{D}}
\def\Aut{\mbox{Aut}}
\def\rre{\mbox{Re}}

\ni
\title[Stability of the Albanese fibration on the Cartwright-Steger surface]
{Stability of the Albanese fibration on the Cartwright-Steger surface}
\author[Vincent Koziarz and Sai-Kee Yeung
]
{Vincent Koziarz and Sai-Kee Yeung
}

\begin{abstract}  
{\it   We verify that the Albanese fibration of the Cartwright-Steger surface is stable, answering a problem left open in [CKY].}
\end{abstract}

\address[Vincent Koziarz]{Univ. Bordeaux, IMB, CNRS, UMR 5251, F-33400 Talence, France}
\email{vincent.koziarz@math.u-bordeaux.fr}
\address[Sai-Kee Yeung ]{Purdue University, West Lafayette, IN 47907, USA
}
\email{yeung@math.purdue.edu
}

\thanks{\noindent{
The second author was partially supported by a grant from the National Science Foundation
}}

\ni{\it }
\maketitle

\bs

The purpose of this note is to resolve a question left open in [CKY] about the semi-stability of the Cartwright-Steger surface, denoted by 
$X$. 
The Cartwright-Steger surface is a smooth complex two ball quotient of Euler number $3$ and first Betti number $b_1=2$.  
It is known that $3$ is the smallest possible number achiveable as  the Euler number of a smooth surface of general type.  The Cartwright-Steger surface
is the only such surface with a non-trivial first Betti number.  The other smooth surfaces of general type with Euler number $3$ are fake projective planes classified in
[PY], [CS1], which have vanishing first Betti numbers.
We refer the readers to
[CS1], [CS2], [CKY] and [BY] for basic geometric properties of $X$.


    Since $h^{1,0}(X)=1$, there is a non-trivial Albanese map.
   Let $\alpha:X\rightarrow E$ be the Albanese fibration, where $E$ is an elliptic curve and is the Albanese
variety of ~$X$.     It is proved in [CKY] that the fibration is reduced and
the genus of a generic fiber is $19$.  A natural question left open in [CKY] is whether the
Albanese fibration is semi-stable, cf. Remark 5.6 of [CKY].   Once the fibration is proved to be semi-stable, it follows from the fact that we are considering a complex two ball quotient with a fibration as above that the fibration is stable.  

The problem turns out to be subtle, and defies conventional algebraic geometric methods
after repeated attempts.  We have to combine the group theoretical  results of [CKY] and explicit equations of [BY] to achieve the purpose.

\begin{theo}
The Albanese mapping $\alpha:X\rightarrow E$ is a stable fibration in the sense of Deligne-Mumford.   
\end{theo}


Recall the following facts about $X$ from [CKY].   The automorphism group is given by $\Aut(X)=\bZ_3$ and the fixed point set of a generator of $\Aut(X)$ consists of $9$ isolated
points, three of type $\frac13(1,1)$ denoted by $O_i, i=1,2,3$, all lying in one fiber of $\alpha$, and six of type $\frac13(1,2)$, denoted by $Q_j, j=1,\dots,6,$
distributed evenly among two different fibers of $\alpha$.  

\begin{lemm} ([CKY])
Assume that $\alpha$ is not stable.  Then there is only one singular fiber, and there is exactly one tacnode singularity on the singular fiber.
The singularity is then one of the points $Q_j, j=1.\dots,6$ listed above.
\end{lemm}

In [BY], explicit equations are found to describe
$X$ as a surface in $P_{\bC}^9$ with coefficients in the rational number field $\bQ$, consisting of $85$ polynomial equations of degree $2$ and $3$.
   As a result, the complex conjugation of $X$ with respect to the natural
complex structure gives a complex surface biholomorphic to itself.  

\begin{lemm} ([BY])
{Denote by $\tau$ the restriction of the complex conjugation of $P_{\bC}^9$ to $X$.  Then $\tau:X\rightarrow X$ is a diffeomorphism of $X$ and the fixed point set is a real subvariety of $X$.}
\end{lemm}

Let $F$ be the fixed point set of $\tau$ on $X$.

\begin{lemm}
The set $F$ contains $O_i, i=1,2,3$, but does not contain $Q_j, j=1,\dots,6$.
\end{lemm}

\ni{\bf Proof}  The action of $\bZ_3$ on $X$ can be described as follows.
According to Remark 5.3 of [BY], the generator $g_3$ acts by
\begin{eqnarray}
&&g_3[U_0:U_1:U_2:U_3:U_4:U_5:U_6:U_7:U_8:U_9]\nonumber\\
&=&[U_0:U_1:U_2:U_3:a U_4:a U_5:a U_6:b U_7:b U_8:b U_9],
\end{eqnarray}
where $a=\zeta_3$ is a primitive cube root of unity, and $b=a^2$.

It follows that the fixed points of $g_3$ have one of the following types,\\
(1): $[U_0:U_1:U_2:U_3:0:0:0:0:0:0]$,\\
(2): $[0:0:0:0:U_4:U_5:U_6:0:0:0]$,\\
(3): $[0:0:0:0:0:0:0:U_7:U_8:U_9]$

From Remark 5.3 of [BY], we already know that the points $O_i$, $i=1,2,3$ are given by the three points $[0:0:0:0:0:0:0:1:0:0], [0:0:0:0:0:0:0:0:1:0], [0:0:0:0:0:0:0:0:0:1]$,
corresponding to points of type (3).  
 In terms of the $85$ explicit equations defining $X$ given in [BY],  the manifold $X$ is the set of complex solutions of the equations
 $q_k(U)=0, k=1,\dots,85$ for $U\in P_{\bC}^9$.  The subvariety $F$ is the set of real solutions
 of the $85$ equations.
 Trying to find solutions of the types (1) to (3), we use Matlab and find that apart of $O_i, i=1,2,3$ as above, there are six other points of type (1), which
 correspond to the points $Q_j, j=1,\dots,6$.  There is no fixed point of type (2).
 The coordinates $U=[U_0:U_1:U_2:U_3:U_4:U_5:U_6:U_7:U_8:U_9]$ of $Q_j, j=1,\dots,6$  are given in terms of  three conjugate pairs of
algebraic numbers, with $U_2/U_0$ given by 
$$-0.0927 \pm 0.1987 \sqrt{-1}, \ \pm 32.0785 \sqrt{-1}, \ 0.0927\pm 0.1987 \sqrt{-1}$$
in decimals respectively.   Hence we may assume that $Q_2=\overline {Q_1}, Q_4=\overline{Q_3}, Q_6=\overline{Q_5}$ after renaming if necessary.
The coordinates were found with the help of Ling Xu using Matlab.
Since the coordinates are not real, the points do not lie on $F$.  

Alternatively, the contradiction in the last step can be reached as follows.  It follows from the information posted on [B] that $U_3^2$ satisfies
the equation $z^3 + (15597z^2)/32 + (641385z)/4096 + 884547/4096 = 0$, from which one easily sees that $U_3$ cannot be real.
\qed

\ms
 Recall that the Albanese map $\alpha:X\rightarrow E$ is defined by $\alpha(\xi)=\int_{x_o}^\xi \omega\pmod \Lambda$, where 
$\omega$ is the holomorphic one form on $X$, $x_o\in X$ is a fixed point on $X$, the integration is taken over any path on $X$ and $\Lambda$ is the lattice generated by $\int_{\gamma}\omega$
as $\gamma$ varies over loops of $X$ generated by $\pi_1(X)$.  The lattice $\Lambda$ takes care of the ambiguity in the choice of path of integration in
the definition above.  Choose $x_o$ to be a fixed point on $X$ with real coordinates, taken to be  $O_j$ for some $j=1,2,3$.

Recall also that the Albanese variety $\alpha(X)=E$ is given by $\bC/(\bZ+\zeta_3\bZ)$, where $\zeta_3$ is a primitive  cube root of unity.  This follows from the discussions
in [CKY].  From the classification of finite automorphism groups of elliptic curves, $E$ has a Weierstrass form given by $z_2^2=z_1^3+1$, cf. [H], page 34 and  [S], 
Theorem III.10.1, which is the unique elliptic curve with automorphism group containing $\bZ_3$.  In this case, $\Aut(E)=\bZ_6=\langle \zeta_6\rangle$, which acts on $E$ by
$\zeta_6\cdot(z_1,z_2)=(\zeta_6^2 z_1, \zeta_6^3 z_2)$.  

In the following, we realize $X$ as a subvariety of $P_{\bC}^{9}$ as given in [BY], and $E$ with the Weierstrass representation in $P_{\bC}^2$ as above.   Then both $X$ and $E$ are invariant by the complex conjugation on the respective projective spaces and the Albanese mapping $\alpha$ will be discussed in terms of the realizations.  We let $U=[U_j]$ and $Z=[Z_\beta]$ be the homogeneous coordinates on the 
ambient manifolds $P_{\bC}^9$ and $P_{\bC}^2$ respectively.  $P_{\bC}^9$ is covered by coordinate charts $\cU_i=\{U|U_i\neq 0\}$ with inhomogeneous coordinates
$\xi_{i,j}=U_j/U_i$ for $0\leqslant j\neq i \leqslant 9$.  Let $\cU_{i,jk}$ be the set of points $p\in X$  so that
the Jacobian matrix of the defining equations of $X$ with respect to $(\xi_{i,j}, \xi_{i,k})$ at $p$ has complex rank $2$.
Then the set of $(\cU_{i,jk}, \xi_{i,j}, \xi_{i,k})$ for $0\leqslant i\leqslant 9$, $0\leqslant j<k\leqslant 9, j\neq i, k\neq i,$ forms a coordinate system on $X$.  For simplicity of notation, we would suppress $i$ and simply denote $\xi_{i,j}$ by $\xi_j$.  The same convention is used for $E$.  Let $\cV_\gamma=\{Z|Z_\gamma \neq 0\}$ and 
$z_{\gamma,\beta}=Z_\beta/Z_\gamma$ for $0\leqslant \beta\neq \gamma \leqslant 2$.  The set of $(\cV_{\gamma,\beta}, z_{\gamma,\beta})$ for $0\leqslant \gamma,\beta \leqslant 2, \beta\neq \gamma$ forms a coordinate system on ~$E$.  Again, we suppress $\gamma$ for notational simplicity.  In affine coordinates
of ~$E$, we may use $z=z_1$ or $z_2$ as coordinate functions in the Weierstrass representation of ~$E$.
Writing $\xi_j=u_j+iv_j$ and $z_\beta=x_\beta+iy_\beta$,  the Albanese map $\alpha$ is given by
$\alpha(\xi)=z(\xi)=x(u,v)+iy(u,v),$ regarding $x, y$ as functions of $u$ and $v$.  Although the construction of these  explicit coordinates is elementary, we explain it in details since we have to use them in the proof below.

\begin{lemm} 
In terms of the realization above, the Albanese map $\alpha$ satisfies the following.\\
(a).  $\bar\alpha(\bar\xi)=\sigma(\alpha(\xi))$ for some $\sigma\in \bZ_6$ in the automorphism group of $E$ fixing the origin.\\
(b).  Let $p\in E$.  Then the complex conjugate of a fiber of $\alpha$ at $p$ is the fiber at $\sigma^{-1}(\bar p)$ for some $\sigma\in \bZ_6$.\\
\end{lemm}

\ni{\bf Proof} 
We describe the Albanese map in terms of the coordinates above.
From the fact that $\alpha$ is holomorphic, we know that $x_{u_j}=y_{v_j}, x_{v_j}=-y_{u_j}$ from Cauchy-Riemann equations, where $j=1,2$.
It follows that $x_{u_j}=(-y)_{(-v_j)}, x_{(-v_j)}=-(-y)_{u_j}$.  This implies that $\bar\alpha(\bar\xi)$ is holomorphic in $\xi$ and hence is a non-trivial
holomorphic map of $X$ to its Albanese torus $E$.   From the functorial properties of the Albanese map, such a holomorphic map is unique up to a biholomorphism of ~$E$.
Hence once the based point $x_o=O_i$ is fixed, $\alpha$ is unique up to $\sigma\in \bZ_6$. 
Hence $\bar\alpha(\bar\xi)=\sigma\circ \alpha(\xi)$.  

(a) leads to $\sigma\circ\alpha(\bar \xi)=\bar\alpha(\xi)$ or $\bar\xi\in\alpha^{-1}(\sigma^{-1}(\bar\alpha(\xi)))$.  Hence by considering the set of $\xi$ with $\alpha(\xi)=p$,
the complex conjugate $\overline{\alpha^{-1}(p)}$ of the fiber at a point $p\in E$ is actually the fiber of $\alpha$ 
at $\sigma^{-1}(\bar p)\in E$.  Hence (b) follows.\
\qed

\ms
\ni{\bf Proof of Theorem 1}  


 Assume on the contrary that the Albanese fibration is not semi-stable, so that according to Lemma 1, there is a tacnode singularity at $Q_j$  on the fiber $L=\alpha^{-1}(p)$ (with 
$p$ one of the fixed points 
of the $\bZ_3$ action on $E$).
From Lemma 4b, $\tau(L)$ is a fiber of $\alpha$ and hence, $\tau(Q_j)$ which is different from $Q_j$ by Lemma 3 is another tacnode singularity, which is a contradiction with Lemma 1.
 \qed


\bigskip

We remark that at the finishing touch of the paper, it is drawn to our attention that there is another proposed independent proof of the stability of the Albanese fibration by Carlos Rito [R], which also relies on the work of [CKY] and [BY], but in a different manner.

\bigskip

\ni {\it Acknowledgement:}   The authors would like to thank Donald Cartwright and Lev Borisov for 
helpful discussions, and to thank Ling Xu for her help with Matlab.  The second author
thanks the hospitality of  Univ. Bordeaux, at which the research here was partly 
conducted.

\bigskip
\noindent{\bf References} 

\ms
\ni [B] Borisov, L., http://www.math.rutgers.edu/$\sim$borisov/CS/

\ms
\ni [BY] Borisov, L., Yeung, S.-K., Explicit equations of the Cartwright-Steger surface, arXiv:1806.08281.

\ms
\ni [CKY] Cartwright, D., Koziarz, V., Yeung, S.-K, On the Cartwright-Steger surface,  J. Algebraic Geom. 26 (2017), 655--689.

\ms
\ni [CS1] Cartwright, D., Steger, T., Enumeration of the $50$ fake projective planes, C. R. Acad. Sci. Paris, Ser. 1, 348 (2010), 11--13,

\ms
\ni [CS2] Cartwright, D., Steger, T., Finding generators and relations for groups acting on the hyperbolic ball, preprint.

\ms
\ni [H]  Husem\"oller, D., Elliptic curves, GTM 106, Springer-Verlag.

\ms
\ni [PY] Prasad, G., Yeung, S-K.,  Fake projective planes, Inv. Math. 168 (2007), 321--370; Addendum, Invent. Math. 182 (2010), 213--227.

\ms
\ni [R] Rito, C., Surfaces with canonical maps of maximal degree, arXiv 1903:03017v1.

\ms
\ni [S] Silverman, J. H., The arithmetics of elliptic curves, GTM 111, Springer-Verlag.


\end{document}